\documentclass{elsart3-1} 
 
\usepackage{epsfig} 
\usepackage{amssymb} 
\usepackage[english,francais]{babel}   
 
\newtheorem{theorem}{Theorem}[section] 
\newtheorem{lemma}[theorem]{Lemma} 
\newtheorem{e-proposition}[theorem]{Proposition} 
\newtheorem{corollary}[theorem]{Corollary} 
\newtheorem{e-definition}[theorem]{Definition\rm} 
\newtheorem{remark}{\it Remark\/} 
 
\newtheorem{theoreme}{Th\'eor\`eme}[section] 
 
\newtheorem{proposition}[theoreme]{Proposition}

\renewcommand{\theequation}{\arabic{equation}} 
\newcommand{\Z}{{\mathbb{Z}}}   
\newcommand{\C}{{\mathbb{C}}}   
\newcommand{\Q}{{\mathbb{Q}}}   
\renewcommand{\P}{{PB(\Sigma, n)}}   
   
\newcommand{\B}{{B(\Sigma, n)}}   
\newcommand{\Bgp}{{B(\Sigma_{g,p}, n)}}   
\newcommand{\U}{\mathcal{U}}   
   
\newcommand{\fd}{{\pi_1(\Sigma)}}

\newcommand{\R}{\mathcal{R}}

\def\og{\leavevmode\raise.3ex\hbox{$\scriptscriptstyle\langle\!\langle$~}} 
\def\fg{\leavevmode\raise.3ex\hbox{~$\!\scriptscriptstyle\,\rangle\!\rangle$}}

\begin{document} 
 
\begin{frontmatter} 
 
\selectlanguage{english} 
\title{Braids on surfaces and finite type invariants} 
 
\vspace{-2.6cm} 
 
\selectlanguage{francais} 
\title{Tresses sur les surfaces et invariants de type fini } 
 
 
\selectlanguage{english} 
\author[a1]{Paolo Bellingeri} 
\ead{bellinge@math.univ-montp2.fr}, 
\author[a2]{Louis Funar} 
\ead{funar@fourier.ujf-grenoble.fr }
\ead[url]{http://www-fourier.ujf-grenoble.fr/$\sim$funar}
\address[a1]{Math\'ematiques, cc 051, 
Univ. Montpellier II, Place Eug\`ene Bataillon, 
34095 Montpellier cedex 5} 
\address[a2]{Institut Fourier, BP 74, Univ.Grenoble I, 
Math\'ematiques,  38402 Saint-Martin-d'H\`eres cedex} 
 
\begin{abstract} 
We prove that there is no functorial universal finite type    
invariant for braids in  $\Sigma\times I$  
if the genus of $\Sigma$ is positive. {\it To cite this article: 
P.Bellingeri, L.Funar, C. R. 
Acad. Sci. Paris, Ser. I 336 (2003).} 
 
\vskip 0.5\baselineskip 
 
\selectlanguage{francais} 
\noindent{\bf R\'esum\'e} 
\vskip 0.5\baselineskip 
\noindent 
Nous d\'emontrons qu'il n'y a pas d'invariant universel    
fonctoriel de type fini pour les tresses dans  $\Sigma \times I$, lorsque    
  $\Sigma$ est une surface orientable de genre positif.  {\it Pour citer cet article~: P.Bellingeri, L.Funar, C. R. 
Acad. Sci. Paris, Ser. I 336 (2003).} 
 
\end{abstract} 
\end{frontmatter} 
 
\selectlanguage{francais} 
\section*{Version fran\c{c}aise abr\'eg\'ee} 
Soit $\Sigma$ une surface compacte, connexe et orientable.  
Le groupe $\B$ de tresses \`a $n$ brins sur $\Sigma$  
est une g\'en\'eralisation naturelle du groupe de tresses d'Artin $B_n$  
et du groupe fondamental de $\Sigma$.  
Ces groupes sont apparus pour la premi\`ere fois  dans l'\'etude des espaces  
de configurations (\cite{FN1}, voir aussi \cite{Bi}), bien que  
certains  \'etaient d\'ej\`a connus avant (le  
groupe $B(S^2, n)$ a \'et\'e introduit par Hurwitz).   
Des pr\'esentations pour ces groupes ont \'et\'e d'abord trouv\'ees  
par Scott  (\cite{Sc}) dans le cas des surfaces ferm\'ees. Celles-ci  
ont \'et\'e  ensuite simplifi\'ees  
par Gonz\'alez-Meneses (\cite{GM}), et am\'elior\'ees et  
g\'en\'eralis\'ees pour les surfaces \`a bord  
dans  \cite{Be}.

Le $\Z$-module  libre engendr\'e par des objets  
$1$-dimensionels plong\'es (comme les tresses, les entrelacs ou  encore 
les enchev\^etrements) admet une filtration naturelle 
obtenue par la r\'esolution des objets singuliers ayant un nombre fini  
de points doubles. L'alg\`ebre gradu\'ee associ\'ee  
(appel\'ee alg\`ebre de diagrammes)  
peut \^etre explicitement calcul\'ee et elle a des  propri\'et\'es 
remarquables de  finitude. Dans ce contexte, un invariant universel est  
une \emph{application}  
de notre cat\'egorie  dans une completion de l'alg\`ebre de 
diagrammes, qui induit un isomorphisme au niveau  
gradu\'e. Le th\'eor\`eme fondamental de la th\'eorie  
de Vassiliev  pour les n\oe uds (mais aussi pour les 
enchev\^etrements, donc les tresses en particulier)  dans 
$\mathbb{R}^3$ est la construction,  
due \`a Kontsevich, d'un invariant universel de type fini \`a  
coefficients dans $\Q$. Un ingredient essentiel est l'existence d'un 
associateur de Drinfel'd \`a coefficients rationnels. 
Une propri\'et\'e remarquable de l'integrale de Kontsevich est  
sa {\em fonctorialit\'e} (\cite{Le-Mu}), permettant par exemple de la  
d\'efinir pour les entrelacs en utilisant l'invariant pour les 
tresses.  
Dans le cas de tresses classiques on sait construire  
un invariant universel de type fini \`a  
coefficients dans $\Z$ (\cite{Pa}), mais on ne sait pas s'il existe  
un tel invariant qui soit aussi  {\em multiplicatif}.  
  
Gonz\'alez-Meneses et Paris (\cite{GMP}) ont construit un invariant universel pour 
les tresses sur les surfaces (ferm\'ees) mais qui n'est pas  
fonctoriel.  
Le propos de cette  note est de d\'emontrer que en effet  
ce r\'esultat ne peut pas \^etre am\'elior\'e:

\begin{theoreme} 
Il n'existe pas d'invariant universel    
fonctoriel de type fini pour les tresses dans  $\Sigma \times I$, lorsque    
  $\Sigma$ est une surface compacte, orientable de genre $g\geq 1$.    
\end{theoreme}   
  
\vspace{0.3mm}  
  
\noindent  
La preuve s'appuie sur la forme explicite des relations  
dans $\B$ et elle  ne d\'epend pas du choix de l'anneau de base 
consid\'er\'e.  
 
\selectlanguage{english} 

\section{Introduction}   
Let $\Sigma$  be a compact, connected and orientable surface.   
The group $\B$ of braids on $n$ strands over $\Sigma$  
is a natural generalization of  both the classical braid group $B_n$  
and the fundamental group $\pi_1(\Sigma)$.   
It  appeared first in the study of   
configuration spaces (\cite{FN1}, see also \cite{Bi}).    
Presentations for $\B$ were derived by Scott (\cite{Sc}),   
further improved  by  Gonz\'alez-Meneses   
for closed surfaces (\cite{GM}) and finally given a very simple  
form in \cite{Be}.    
In the case of holed spheres the latter have been   
previously obtained by Lambropoulou (\cite{La}).

Let us consider a category of embedded 1-dimensional objects like   
braids, links, tangles etc. There is a natural filtration on the free   
$\Z$-module generated by the objects, coming from the   
singular objects with a given number of double points.  
The main feature of this filtration is that the   
associated grading, which is called the diagrams algebra,   
can be explicitly computed, and has some salient finiteness properties.    
By  universal finite type invariant one generally means  
a {\em map} from our category   
into some completion of the diagrams algebra,  
which induces an isomorphism at graded level.   
For instance the celebrated Kontsevich integral is such a universal   
invariant. A key ingredient is the existence of a Drinfel'd associator 
with rational coefficients. Notice that there exists a universal 
invariant for usual braids over $\Z$ (\cite{Pa}), but it is not known whether 
there exists a multiplicative one over $\Z$.  
In fact an essential feature of the Kontsevich integral   
is its  {\em functoriality} (\cite{Le-Mu}): it is precisely this property   
which enables one to extend the Chen iterated integrals  
from braids to links (\cite{Le-Mu}).

Gonz\'alez-Meneses and Paris  
constructed a universal invariant for braids on surfaces  
(\cite{GMP}), but their invariant  
is not functorial (i.e. multiplicative).  
The purpose of this note is to show that actually their result   
cannot be improved:

\begin{theorem} \label{main}   
There does not exist any functorial universal  finite type invariant    
for braids in $\Sigma\times I$, if the  surface $\Sigma$ is of  
genus $g\geq 1$.    
\end{theorem}

\begin{remark} 
In particular there does not exist any universal 
invariant for tangles in $\Sigma\times I$  which is  
functorial with respect to the vertical composition 
of tangles. Thus the  Andersen-Mattes-Reshetikhin invariant   
(\cite{AMR}) cannot be made tangle functorial.   
Notice that the functoriality is essential if one seeks for an extension    
to links in arbitrary 3-manifolds. 
The non-existence holds true a fortiori for the category of  
tangles in 3-manifolds with boundary, unless these   
are cylinders over planar surfaces. 
\end{remark} 
 
\section{Preliminaries}   
\noindent {\bf 2.1. Surface braids.} Set  $\Sigma_{g,p}$ for    
the compact orientable surface of genus $g$ with $p$ boundary  
components. We denote by $\sigma_1,...,\sigma_{n-1}$ the standard  
generators of the braid group on a disk embedded in $\Sigma_{g,p}$,   
viewed as elements of $\Bgp$.  Let also   
$a_1,...,a_{g}, b_1,...,b_g, z_1,...,z_{p-1}$ be the generators   
of $\pi_1(\Sigma_{g,p})$, where $z_i$ denotes a loop around the  
$i$-th boundary component.  
Assume that the base point of the fundamental group is   
the startpoint of the first strand. Then each $\gamma\in \{a_1,...,a_{g}, b_1,...,b_g, z_1,...,z_{p-1}\}$ can be realized  by an element,  
denoted also by $\gamma$, in  
$\Bgp$, by considering the braid whose first strand is describing the  
curve $\gamma$ and whose other strands are constant. We denote $[a,b]=a b a^{-1} b^{-1}$.   Following \cite{Be} we have:   
\begin{theorem}  
A presentation for $\Bgp$ ($n\geq 2$) is given by:    
  
$\rm 1.)$ Generators: $\sigma_1,...,\sigma_{n-1}, a_1,...,a_{g}, b_1,...,b_g, z_1, ...,$ $z_{p-1}$.     
   
{\rm 2.i)} Braid relations:   
\begin{center}   
 $ \sigma_{i}\sigma_{i+1}\sigma_i=\sigma_{i+1}\sigma_{i}\sigma_{i+1} \; (i=1, \dots, n-2)  \quad   
 \sigma_i\sigma_j=\sigma_j\sigma_i \;  (\mbox{if } \, \mid i -j\mid \geq   
2). $   
\end{center}   
   
{\rm 2.ii)}   Commutativity relations:   
\begin{center}   
  $ \quad [ a_r, \sigma_i ]= [ b_r, \sigma_i ]=[z_k, \sigma_i]= 1,  \,
  \,  
(i > 1, \, 1\leq r\leq g, \, 1\leq k\leq p-1)\, ;$ 
 
$ [a_r, \sigma_1^{-1} a_r \sigma_1^{-1}] = [b_r, \sigma_1^{-1} b_r \sigma_1^{-1}]= 
[z_k, \sigma_1^{-1} z_k \sigma_1^{-1}]= 1, \,\, (1\leq r\leq g, \,
1\leq k\leq p-1) ;$ 
   
$  [a_r, \sigma_{1}^{-1}a_s\sigma_{1}]= [a_r, \sigma_{1}^{-1}b_s\sigma_{1}]= 
[b_r, \sigma_{1}^{-1} a_s\sigma_{1}]= 
[b_r, \sigma_{1}^{-1} b_s\sigma_{1}]= [z_i, \sigma_{1}^{-1}z_j\sigma_{1}]= 1,$ 
 
$(1\le s < r \le g, \,1\le j<i \le p-1)\, ;$ 
 
$[a_r, \sigma_{1}^{-1}z_k\sigma_{1}]=[b_r,
\sigma_{1}^{-1}z_k\sigma_{1}]=1,   \,\,  (1\leq r\leq g, \,
1\leq k\leq p-1)\, ; $  
 
\end{center}   
 
{\rm 2.iii)} Skew commutativity relations on each handle (when $g>0$): 
  
\begin{center}    
$[ a_r, \sigma_{1}^{-1}b_{r}\sigma_{1}^{-1}]=   
\sigma_{1}^2 , \,\, (1\leq r\leq g)\, ;$   
\end{center}   
   
{\rm 2.iv)}  When $p=0$ we have the additional relation:   
\begin{center}   
$[a_1,b_1^{-1}][a_2, b_{2}^{-1}]...[a_{g},b_{g}^{-1}]=   
\sigma_1\sigma_2...\sigma_{n-2}\sigma_{n-1}^2   
\sigma_{n-2}...\sigma_2\sigma_1\,.    
   $   
\end{center}   
\end{theorem}   
   

\vspace{0.2cm}   
\noindent {\bf 2.2. The Vassiliev filtration.}    
Singular   braids are obtained from braids by allowing  a finite   
number of transverse double points.   
One can associate to each singular surface braid a linear combination  
of braids  by desingularizing each crossing as follows:    
$$   \left(  \raisebox{-4mm}{\psfig{figure=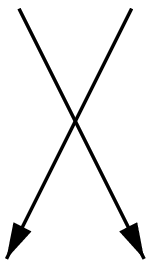,width=15pt}} \right)   \;   
\rightarrow  \; \,\left( \, \raisebox{-4mm}{\psfig{figure=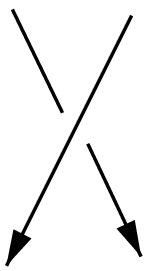,width=15pt}}  \right)     
\; - \; \left( \, \raisebox{-4mm}{\psfig{figure=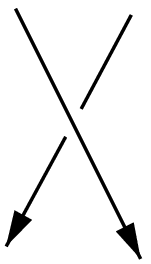,width=15pt}}  \right)  \,.   
$$    
We denote by $\mathcal{V}^d$  the submodule of $ \Z [\B]$   
generated by the desingularizations of (singular) surface braids with   
$d$ double points. Then $\mathcal{V}^d=J^d$, where 
 $J\subset \Z [\B]$ is the two-sided ideal generated   
by $\{ \sigma_i - \sigma_i^{-1}\}_{i=1, \dots, n-1}$.
Set $gr^*\Z[\B]$ for the   
associated graded algebra. 
A morphism of $\Z$-modules 
$u:\Z[\B]\to A$ is a Vassiliev invariant of degree 
$\leq d$ if $u$ vanishes on $\mathcal{V}^{d+1}$.

\begin{e-definition}   
A universal  Vassiliev invariant for braids is a linear map 
$Z: \C[\B] \to {\mathcal A}$ such that, for any finite type 
invariant  $v: \C[\B]\to A$ there exists a homomorphism 
$u_v:{\mathcal A}\to A$ satisfying $u_v\circ Z=v$. The invariant is
functorial if ${\mathcal A}$ has an algebra structure and $Z$ is 
an algebra homomorphism.   
\end{e-definition}

\vspace{0.2cm} 
\noindent {\bf 2.3. Presentations of the algebra $gr^*\Z[\B]$.}
Let $S_n$ denote the symmetric group on $n$ elements. 
The proof given in \cite{GMP} for $p=0$ extends without essential modifications
to the general case: 
\begin{proposition}
The graded algebra $gr^*\Z[\B]$ is the quotient of 
$\Z[Z_{ij}^{\gamma}]_{1\leq i,j\leq n; \gamma\in \fd}\ltimes S_n$,
the semi-direct product of a free non-commutative algebra in the given
generators  with $S_n$,      
by the following ideal of relations $\R$:    
$$ Z_{ij}^{\gamma}=Z_{ji}^{\gamma^{-1}}, \;\,  [Z_{ij}^{\gamma}, Z_{kl}^{\delta}]=0, \,\;  
[Z_{ij}^{\gamma}, Z_{jk}^{\delta}+Z_{ik}^{\gamma\delta}]=0. $$  
\end{proposition}

\noindent This algebra can be described in terms of chord 
diagrams with beads (\cite{GK}). 
Consider the free $\Z$-module generated by the horizontal chord   
diagrams with endpoints on $n$ vertical segments, where the   
vertical intervals between   
two successive chords endpoints are labeled by elements  
of $\Z[\pi_1(\Sigma)]$  
called beads. We obtain a $\Z[\pi_1(\Sigma)]$-algebra ${\mathcal  
  D}(n,\Sigma)$  
by imposing the (4T) relations    
along with (P) below, permitting to   
push beads across a chord:    
   
\vspace{0.3cm}   
$$   
\hspace{-10pt} (4T)\,:  
\hspace{2pt}\raisebox{-8mm}{\psfig{figure=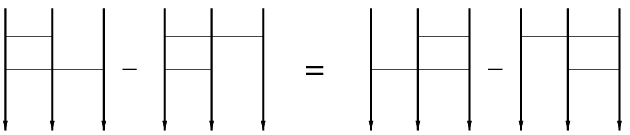,width=7cm}}\;,\;\;\,\;  
(P)\,:\hspace{2pt} \raisebox{-8mm}{\psfig{figure=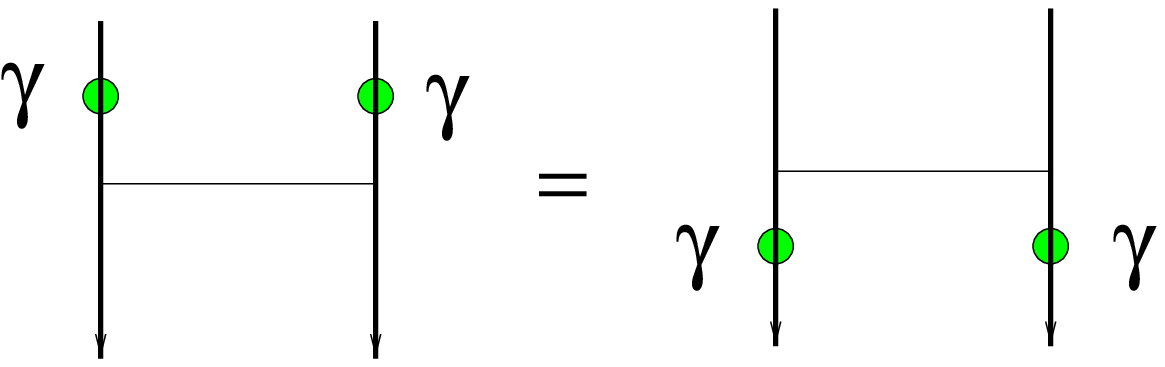,width=4.5cm}}     
\vspace{0.3cm}   
$$   
\noindent Set $F(\Sigma_{g,p})=\{A_{s},B_s, {A_s}^{-1},   
{B_s}^{-1},Z_{\alpha},{Z_{\alpha}}^{-1}\}_{
1\leq s\leq g; 1\leq \alpha\leq p-1}$.
The   
algebra ${\mathcal D}(n,\Sigma)$ is the quotient of the 
free non-commutative algebra   
$\Z[Z_{ij}, \gamma^i]_{1\leq i,j\leq  
  n; \gamma\in F(\Sigma_{g,p})}$ by the ideal of relations:  
$$ [\gamma^i, \delta^j]=0,  \; 
 [\gamma^i + \gamma^j, Z_{ij}]=0,\; \mbox{ for } i\neq j, \;\mbox{ and }
\gamma,\delta\in F(\Sigma_{g,p}), $$  
$$ [\gamma^i, Z_{jk}]=0,  \;  
Z_{ij}=Z_{ji}, \;   
[Z_{ij}, Z_{kl}]=0,  \; 
[Z_{ij}, Z_{jk}+Z_{ik}]=0, \; \mbox{ if } i,j, k  \mbox{ are distinct}. $$  
If $\Sigma$ is closed one adds also the relations:   
$\sum_{s=1}^g [A_{s}^i, B_s^i]=0.$  
\begin{proposition}
The map $gr^*\Z[\B])\to {\mathcal D}(n,\Sigma)\ltimes S_n$   
sending   $Z_{ij}^{\gamma}$   
to $\gamma^i Z_{ij}{\gamma^i}^{-1}$, where  $\gamma^i$ denotes the bead   
labeled $\gamma$ on the $i$-th strand, is an isomorphism.
\end{proposition}    
\begin{remark}
${\mathcal D}(n,\Sigma)$ is isomorphic to the  
algebra of trivalent graphs with beads on oriented edges from \cite{GK}. 
\end{remark}
   
We denote by $H_1(A)$ the abelianization of the graded algebra $A$,   
namely the quotient by the homogeneous ideal generated by the   
$ab-ba$, for $a,b\in A$.

\begin{corollary}\label{coro}
$H_1(gr^*\Z[\B])\cong
\Z[Z_{12}, A^i_1, {A_1^{i}}^{-1}, B^j_1, {B_1^{j}}^{-1},Z_{1}^i,
{Z_{1}^i}^{-1}, \tau]/(\tau^2=1)$ as graded commutative algebras. 
Here $\tau$ corresponds to the signature on $S_n$.   
The degree of $Z_{12}$ is one, while the other generators are of degree zero. 
\end{corollary}

\section{Proof of the Theorem}
\noindent Assume that such a functorial universal invariant $Z$ exists. 
The algebra $M=Z(\C[\B])$ is filtered by the ideals
$\mathcal{V}^jM=Z(\mathcal{V}^j)$ and we set $gr^*(M)$ 
for the associated graded algebra. 
\begin{lemma}
$Z$ induces an isomorphism 
of graded algebras $Z^*:gr^*\C[\B]\to gr^*M$. 
\end{lemma}
\noindent {\em Proof.}\, It suffices to prove that $Z^*$ is
injective. One introduces the pairing \\ 
$\langle,\rangle :gr^dM \times {\rm Hom}(gr^d\C[\B],\C) \to \C$, defined as 
$\langle x,\lambda\rangle=u_{I(\lambda)}(\hat{x})$, where 
$\hat{x}\in \mathcal{V}^dM$ is a lift of $x$ and 
$I(\lambda):\C[\B]\to \C$
is an extension of $\lambda$ as a Vassiliev invariant 
of degree $\leq d$. 

Notice that any $\Z$-modules homomorphism 
$\lambda:\frac{\mathcal{V}^{d}}{\mathcal{V}^{d+1}}\to
A$ can be lifted (non-uniquely) to a Vassiliev invariant 
$I(\lambda)$ of degree $\leq d$, which coincides with $\lambda$ on 
$\frac{\mathcal{V}^{d}}{\mathcal{V}^{d+1}}$. 

If $a\in gr^d\C[\B]\cap  \ker Z^*$, then $\langle Z^*(a),\lambda\rangle=0$ 
for all $\lambda\in {\rm Hom}(gr^d\C[\B],\C)$. Choose a representative 
$\hat{a}\in {\mathcal V}^d$ for $a$. Then 
$\langle Z^*(a),\lambda\rangle=u_{I(\lambda)}(Z(\hat{a}))=I(\lambda)(\hat{a})=
\lambda(a)$. This implies $a=0$ because $gr^d\Z[\B]$ 
is torsion-free (see \cite{GMP}) and our claim follows. $\; \Box$

\begin{lemma}
The map $Z^*_{ab}: H_1(gr^*\C[\B])\to H_1(gr^*M)$ 
induced by $Z^*$ is trivial in degree one.   
\end{lemma}
\noindent {\em Proof.}
The image of $\sigma_1^2-1$ in $\frac{\mathcal{V}^1}{\mathcal{V}^2}$
is the element $Z_{12}$, associated to the 
diagram \raisebox{-5mm}{\psfig{figure=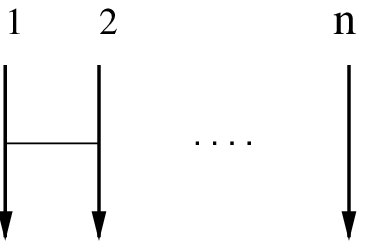,width=50pt}}. 
Furthermore its image under $Z^*$ is 
$Z(\sigma_1^2)-1\in gr^1(M)$. But 
$Z(\sigma_1^2)=Z([a_1,\sigma_{1}^{-1}b_1\sigma_{1}^{-1}])=
1 \in   H_1(gr^*(M))$, since $Z$ is multiplicative and 
the image of a commutator in the abelianization is trivial. 
Therefore $Z^*_{ab}(Z_{12})=0$ and the corollary \ref{coro} 
implies our claim.
$\;  \Box$

The first lemma implies that $Z^*_{ab}$ should be an isomorphism 
$H_1(gr^*\C[\B])\to H_1(gr^*M)$, while the last lemma
shows that this map is trivial in degree one.
But the degree one component of $H_1(gr^*\C[\B])$ is nonzero 
from the  the corollary \ref{coro}, and thus we get a contradiction.  
This settles  the  
theorem.

\section{Comments}  
  
\noindent Let $\P$ be the pure  surface braid group on $n$ strands and    
$K(\Sigma,n)$ the kernel of the   
natural homomorphism $\theta: \P \to \fd^n$ which forgets   
about the braiding and keeps only the  fundamental group   
information of each strand. Then  $J$ coincides with  
the two-sided ideal in $\Z [\B]$ generated by the augmentation ideal  
$I_{K(\Sigma,n)}$
One cannot use the Chen iterated integrals for the subgroup  
$K(\Sigma,n)$  since this group is not of finite type (see 
also \cite{CKX}).   
An interesting alternative would be to replace the Vassiliev  
filtration ${\mathcal V}^*$ by that associated to the   
augmentation ideal $I_{\P}$. The associated graded algebra 
$gr^*_{I_{\P}}\C[\B]$ is now isomorphic to  $\U\P\ltimes S_n$, where 
$\U G$ denotes the universal enveloping algebra 
of the Lie algebra of the group $G$ over $\C$
(see \cite{H}). In this setting one can find a functorial  
universal-type invariant 
 (for the new filtration) $\C[\B]\to \overline{\U\P}\ltimes S_n$  using the  
Chen iterated integrals, which takes values in the completion of 
an algebra of  {\em symplectic chord diagrams}.  
Specifically $\U\P$ is a chord  
diagram algebra as defined above, but now the actions of    
beads on two distinct segments are not anymore assumed to commute,  
and they are subjected to the following relation:   
$$
\raisebox{-8mm}{\psfig{figure=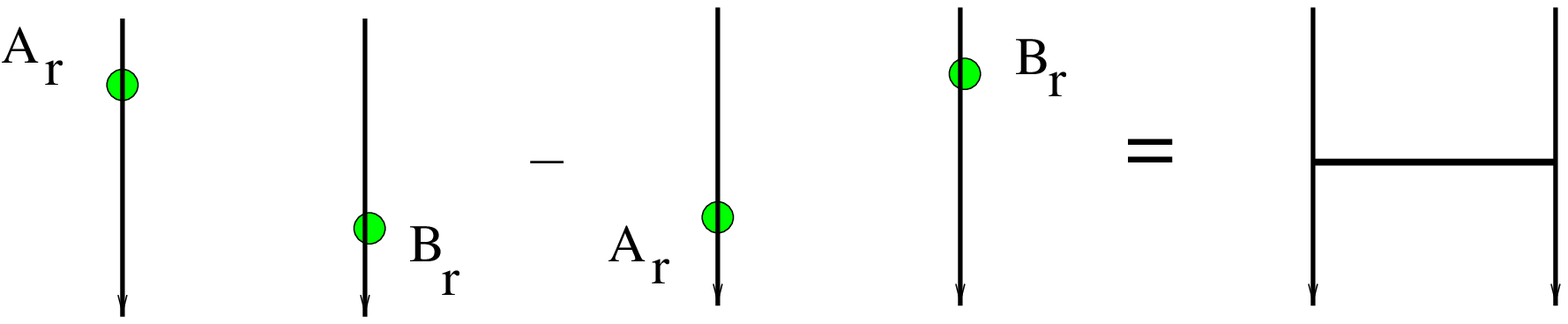,width=7cm}}  
$$    
Nevertheless the grading is different from the previous algebra of 
chord diagrams since beads generators $A_s^i, B_t^j$ are now given the  
degree 1 while each chord has degree 2.    
Alternatively one can present this algebra as  
the quotient of  $\C[A_s^k, B_t^k, Z_{\alpha i}, Z_{ij}]_{  
1\leq s,t \leq g, 1\leq \alpha\leq p, 1\leq i,j,k \leq  
  n}$, where $\deg\;A_s^k=\deg\; B_t^k =1$ and  $\deg \; Z_{mj} =2$  
(unless the case when $g=0$ and one can renormalize the degree   
of the generators $A_{ij}$ to be 1) by the   
following relations:  
  
$\bullet$ The extended infinitesimal braid relations:   
\[ Z_{ij}=Z_{ji},   
[Z_{ij}, Z_{kl}]=0, \mbox{ if } \{i,j\}\cap \{k, l\}=\emptyset,   
[Z_{ij}, Z_{jk}+Z_{ik}]=0,\mbox{ if } i, j, k \mbox{ are  
  distinct}. \]  
\[ [Z_{\alpha j}, Z_{kl}]=0, \mbox{ if } j\not\in \{k, l\},   
 [Z_{\alpha j}, Z_{\beta k}]=0, \mbox{ if } \{\alpha,j\}\cap \{  
\beta, k\}=\emptyset,   
[Z_{\alpha j}, Z_{\alpha k}+Z_{jk}]=0,\mbox{ if } j\neq k. \]  
  
$\bullet$ The relations coming from the fundamental group of $\Sigma$:   
\[ [A_s^{i}, A_{r}^k]=[B_s^{i}, B_{r}^k]=0, \mbox{ if } i\neq k,   
[A_s^i, B_r^j]=0, \mbox{  if } r\neq s \mbox{ and } i\neq j, \]  
\[ \sum_{s=1}^g[A_s^k, B_s^k] + \sum_{j=1}^n Z_{jk} +  
\sum_{\alpha=n+1}^{n+p} Z_{\alpha k} =0. \]  
  
$\bullet$ The mixed relations:   
\[ [Z_{jk}, A_s^i]=0, \mbox{ if } i\not\in \{j,k\},   
 [Z_{\alpha k}, A_s^i]=0, \mbox{ if } i\neq k,   
 [A_s^j + A_s^k, Z_{jk}]=0,   
 [B_s^j + B_s^k, Z_{jk}]=0. \]  
  
$\bullet$ The twist relation (making sense only for $g\geq 1$):  
 \[ [A_s^i, B_s^j]=Z_{ij}, \mbox{  if } i\neq j. \]  
The closed case is proved in \cite{Bez}, for the general 
case see (\cite{H}, Thm.12.6).
A different approach based on the weight filtration   
on $\P$ is given in \cite{NT}.   
It is not difficult to see that the completion of this symplectic   
chord diagrams   
algebra surjects onto $\overline{\U K(\Sigma,n)}\otimes  
\overline{\U(\pi_1(\Sigma))}$   
and hence it furnishes an invariant in the usual   
chord diagrams but whose coefficients are now   
formal series from $\overline{\U(\pi_1(\Sigma))}$ instead of merely  
elements of $\pi_1(\Sigma)$. The relation of the  
latter with the Vassiliev filtration seems quite obscure.   
  
\noindent It is worth mentioning that the 
universal flat Chen connection whose monodromy yields this invariant is 
{\em not quadratic} since the respective configuration spaces are  
not formal unless $g=0$ or $n=2$ (see \cite{Bez}). 
  
\begin{remark}   
In the case of surface pure braid groups   
the exact sequence  
$$   
\quad 1 \to \pi_1(\Sigma\setminus  \{ (n-1) \; \mbox{points} \}) \to \P \to PB(\Sigma,n-1)  \to  1 \,    
$$   
is not split unless $\Sigma$ is a torus or    
$\Sigma$ has non-empty boundary (see \cite{GG}).    
Moreover even when  split, the groups $PB(\Sigma,n)$ are not   
{almost-direct products} if $g\geq 1, n\geq 3$.     
This explains why the arguments in \cite{KO}   
fail for higher genus, as already noticed in    
\cite{Bez}. In particular it is still unknown whether    
$gr^*PB(\Sigma,n)$ is torsion-free or if $\P$ is residually nilpotent,  
which would shed some light about the relevance of the filtration 
considered in this section.  
\end{remark}

\section*{Acknowledgements} 
The authors are indebted to J.Birman, M.Eisermann, J.Gonz\'alez-Meneses,  
J.Guaschi, T.Kohno, T.Le, S.Lambropoulou,  D.Matei, L.Paris and    
S.Papadima for useful discussions and to the referee 
for suggesting us the present  definition 
of a functorial universal  Vassiliev invariant.


\begin{thebibliography}{00} 
\selectlanguage{english}   
  
  
\bibitem{AMR}   
J.E.~Andersen, J.~Mattes and N.~Reshetikhin,  
\newblock Quantization of the algebra of chord diagrams,   
\newblock {\em Math. Proc. Camb. Philos. Soc.} 124 (1998), 451-467.   
   
   
\bibitem{Bez}   
R.~Bezrukavnikov, {Koszul DG-algebras arising from configuration   
  spaces}, {\em G.A.F.A.} 4 (1994), no. 2, 119--135.   
   
   
   
   
   
\bibitem{Be}   
P.~Bellingeri,  {On presentation of surface braid groups},    
math.GT/0110129.   
   
   
   
\bibitem{Bi}   
J.~Birman, 
\newblock  {\em Braids, links, and mapping class groups},   
\newblock Ann.Math.Studies 82, Princeton, 1973.   
   
   
   
  
   
   
   
   
%
   
   

\bibitem{CKX} 
F.R.~Cohen, T.~Kohno and M.A.~Xicotencatl, 
{Orbit configuration spaces associated to discrete subgroups of
  PSL(2,R)}, math.AT/0310393. 



   
\bibitem{FN1}   
E.~Fadell and L.~Neuwirth,
Configuration spaces, {\em Math.Scandinavica}, 10(1962), 111-118.    
   
   
   
   
   
\bibitem{GK}   
S.~Garoufalidis, A.~Kricker,   
{A rational noncommutative invariant of boundary links}, math.GT/0105028.    
   
   
\bibitem{GM}   
J.~Gonz\'alez-Meneses, New presentations of surface braid groups,   
{\em J.K.T.R.} 10 (2001), no. 3, 431--451.    
   
   
   
\bibitem{GMP}   
J.~Gonz\'alez-Meneses and L.~Paris,    
\newblock {Vassiliev invariants for braids on surfaces},    
\newblock {\em T. A. M. S.} 356(2004), 219-243.    
   
\bibitem{GG}   
D.L.~Goncalves and J.~Guaschi,    
{On the structure of surface pure braid groups},  
{\em J. Pure Appl. Algebra} 182 (2003),  33--64.    
   
\bibitem{H}   
R.~Hain, {Infinitesimal presentations of the Torelli groups},
{\em J.A.M.S.} {10} (1997), 
  597--651.

   
   
   
   
\bibitem{KO}   
T.~Kohno and T.~Oda,     
{The lower central series of the pure braid group of an algebraic curve},    
201--219,   {\em Adv. Stud. Pure Math.}, 12 (1987), 201--219.    
   
   
   
   
   
  
   
\bibitem{La}   
S.~Lambropoulou, {Braid structures in knot complements, handlebodies  
  and 3-manifolds}, {\em Knots in Hellas '98 (Delphi)}, 274--289, Ser.   
Knots Everything, 24, World Sci. Publishing, River Edge, NJ, 2000.    
   
   
\bibitem{Le-Mu}   
T.Q.T.~Le and J.~Murakami.   
\newblock Representations of the category of tangles by {K}ontsevich's iterated integrals.   
\newblock {\em Commun.Math.Phys.}, 168 (1995), 535--563.   
   
  
\bibitem{NT}  
H.~Nakamura and N.~Takao, {Galois rigidity of pro-$l$ pure braid  
  groups of algebraic curves}, {\em T.A.M.S.} 350(1998), 1079-1102.   
  
   
\bibitem{Pa}   
S.~Papadima,    
The universal finite-type invariant for braids, with integer    
coefficients, {\em Topology Appl.}, 118 (2002), 169-185.    
   
   
   
\bibitem{Sc}   
G.P.~Scott.   
\newblock  Braid groups and the group of homeomorphisms of a surface.   
\newblock {\em Math.Proc.Cambridge Phil. Soc.}, 68 (1970), 605--617.   
   
  

 
\end{thebibliography}
\end{document}